\documentclass[12pt,leqno]{article}
\usepackage{amssymb}
\usepackage[mathscr]{eucal}
\usepackage{amsmath,amssymb,latexsym,theorem,bbm}
\newcommand{\EE}{\mathsf{E}}

\newcommand{\NN}{\mathbb{N}}
\newcommand{\PP}{\mathsf{P}}
\newcommand{\QQ}{\mathbb{Q}}
\newcommand{\RR}{\mathbb{R}}

\newcommand{\cA}{{\mathcal A}}
\newcommand{\cB}{{\mathcal B}}

\newcommand{\cG}{{\mathcal G}}

\newcommand{\dd}{\mathrm{d}}
\newcommand{\ee}{\mathrm{e}}

\DeclareMathOperator*{\argmax}{arg\,max}

\newcommand{\cov}{\operatorname{Cov}}

\newcommand{\halpha}{\widehat{\alpha}}

\renewcommand{\leq}{\leqslant}
\renewcommand{\geq}{\geqslant}

\newcommand{\proofend}{\hfill\mbox{$\Box$}}

\numberwithin{equation}{section}

\theoremstyle{change} \theorembodyfont{\em}
\newtheorem{Lem}{Lemma.}[section]
\newtheorem{Thm}[Lem]{Theorem.}

\theorembodyfont{\rm}
\newtheorem{Rem}[Lem]{Remark.}

\begin{document}

\begin{center}
 {\bfseries\Large $\alpha$-Wiener bridges: singularity of induced measures \\[2mm]
                  and sample path properties}\\[5mm]

 {\sc\large M\'aty\'as $\text{Barczy}^{*,\diamond}$} {\large and}
 {\sc\large Gyula $\text{Pap}^*$}
\end{center}

\vskip0.2cm

* University of Debrecen, Faculty of Informatics, Pf.~12, H--4010 Debrecen, Hungary;
 e--mail: barczy@inf.unideb.hu (M. Barczy), papgy@inf.unideb.hu (G. Pap).

$\diamond$ Corresponding author.

{Running head: $\alpha$-Wiener bridges: singularity of induced measures}



\renewcommand{\thefootnote}{}

\vspace*{-5mm}

\begin{abstract}
Let us consider the process \ $(X_t^{(\alpha)})_{t\in[0,T)}$ \ given by the SDE
 \ $\dd X_t^{(\alpha)}=-\frac{\alpha}{T-t}\,X_t^{(\alpha)}\,\dd t+\dd B_t$, \ $t\in[0,T)$,
 \ where \ $\alpha\in\RR$, \ $T\in(0,\infty)$, \ and \ $(B_t)_{t\geq 0}$ \ is a standard
 Wiener process.
In case of \ $\alpha>0$ \ the process \ $X^{(\alpha)}$ \ is known as an \ $\alpha$-Wiener bridge,
 in case of \ $\alpha=1$ \ as the usual Wiener bridge.
We prove that for all \ $\alpha, \beta\in\RR$, \ $\alpha\ne\beta$, \ the probability measures
 induced by the processes \ $X^{(\alpha)}$ \ and \ $X^{(\beta)}$ \ are singular on
 \ $\big(C[0,T),\cB(C[0,T))\big)$.
\ Further, we investigate regularity properties of \ $X_t^{(\alpha)}$ \ as \ $t\uparrow T$.
\end{abstract}


\section{Introduction}

There has been a lot of work concerning questions of absolute continuity and singularity for
 various types of stochastic processes on finite and infinite time intervals, see for example
 Jacod and Shiryaev \cite{JacShi}, Ben-Ari and Pinsky \cite{BenPin}, and Prakasa Rao \cite{Rao_2008}.
However, most of the literature deal with time homogeneous diffusion processes.
In this paper we study absolute continuity and singularity of $\alpha$-Wiener bridges which are
 time inhomogeneous diffusion processes.
We also present some results on the sample path behavior of these processes.

Let \ $T\in(0,\infty)$ \ be fixed.
For all \ $\alpha\in\RR$, \ we consider the process \ $(X_t^{(\alpha)})_{t\in[0,T)}$ \ given
 by the stochastic differential equation (SDE)
 \begin{align}\label{alpha_W_bridge}
  \begin{cases}
   \dd X_t^{(\alpha)}=-\frac{\alpha}{T-t}\,X_t^{(\alpha)}\,\dd t+\dd B_t,\qquad t\in[0,T),\\
   \phantom{\dd} X_0^{(\alpha)}=0,
  \end{cases}
 \end{align}
 where \ $(B_t)_{t\geq 0}$ \ is a $1$-dimensional standard Wiener process
 on a probability space \ $(\Omega,\cA,\PP)$.
\ To our knowledge, these kind of processes have  been first considered by
 Brennan and Schwartz \cite{BreSch}, and see also Mansuy \cite{Man}.
In Brennan and Schwartz \cite{BreSch} the SDE \eqref{alpha_W_bridge} is used to model
 the arbitrage profit associated with a given futures contract in the absence of transaction costs.
In case of \ $\alpha>0$ \ the process \ $X^{(\alpha)}$ \ is known as an $\alpha$-Wiener bridge,
 in case of \ $\alpha=1$ \ as the usual Wiener bridge.
By formula (5.6.6) in Karatzas and Shreve \cite{KarShr}, the SDE \eqref{alpha_W_bridge}
 has a unique strong solution, namely,
 \begin{align}\label{alpha_W_bridge_intrep}
   X_t^{(\alpha)}=\int_{0}^t\left(\frac{T-t}{T-s}\right)^\alpha\,\dd B_s,\qquad t\in[0,T),
 \end{align}
  defined on a filtered probability space
 \ $\big(\Omega, \cA, (\cA_t)_{t\in[0,T)}, \PP\big)$ \ constructed by the help of
 the standard Wiener process \ $B$, \ see, e.g., Karatzas and Shreve \cite[Section 2.5.A]{KarShr}.
This filtered probability space satisfies the so called usual conditions, i.e.,
 \ $(\Omega,\cA,\PP)$ \ is complete, the filtration \ $(\cA_t)_{t\in[0,T)}$ \ is right-continuous,
 \ $\cA_0$ \ contains all the $\PP$-null sets in \ $\cA$ \ and \ $\cA=\cA_{T-},$ \ where
 \ $\cA_{T-}:=\sigma\left(\bigcup_{t\in[0,T)}\cA_t\right)$.

In Section \ref{Section_Preliminaries} \ we calculate the covariance of \ $X^{(\alpha)}_s$ \ and
 \ $X^{(\beta)}_t$ \ for all \ $s,t\in[0,T)$ \ and \ $\alpha$, $\beta\in\RR$.
\ Further, we recall a strong law of large numbers and a law of the iterated logarithm for
 continuous local martingales which will be used for proving regularity properties of \ $X^{(\alpha)}$.

In Section \ref{Section_sample_paths} we prove that \ $X_t^{(\alpha)}\to 0$ \ almost surely as \ $t\uparrow T$
 \ in case of \ $\alpha>0$, \ see, Lemma \ref{LEMMA5}, and that's why we can use the expression
 '$\alpha$-Wiener bridge' for \ $X^{(\alpha)}$ \ in this case.
Lemma \ref{LEMMA5} can be considered as a generalization of Lemma 5.6.9 in Karatzas and Shreve \cite{KarShr}.
We will also examine what happens in case of \ $\alpha\leq 0$, \ see Remark \ref{REMARK2}.
Further, we investigate regularity properties of \ $X_t^{(\alpha)}$ \ as \ $t\uparrow T$.
\ In case of \ $\alpha\geq \frac{1}{2}$ \ we have theorems of type of the law of the iterated logarithm
 for \ $X^{(\alpha)}$, \ see Theorems \ref{THM_path1} and \ref{THM_path2}.

In Section \ref{Section_measures} we investigate the absolute continuity and singularity
 of the probability measures induced by the processes \ $X^{(\alpha)}$ \ with different values of \ $\alpha$.
\ Namely, we show that for all \ $\alpha, \beta\in\RR$, \ $\alpha\ne\beta$, \
 the probability measures induced by the processes \ $X^{(\alpha)}$ \ and \ $X^{(\beta)}$
 \ on \ $\big(C[0,T),\cB(C[0,T))\big)$ \ are singular, where \ $C[0,T)$ \ is the space of
 continuous functions from \ $[0,T)$ \ into \ $\RR$, \ and \ $\cB(C[0,T))$ \ denotes the
 Borel $\sigma$-algebra on \ $C[0,T)$, \ see Theorem \ref{THM_singular2}.
We note that Prakasa Rao \cite[Theorem 5]{Rao_2008} proved a similar statement for fractional Wiener processes.
Namely, he showed that if \ $(W_{H_i}(t))_{t\geq 0}$, \ $i=1,2$, \ are two fractional Wiener processes with
 Hurst indices \ $H_1,H_2\in(0,1)$, \ $H_1\ne H_2$, \ then the probability measures induced by
 the processes \ $W_{H_i}$, \ $i=1,2$, \ on \ $\big(C[0,\infty),\cB(C[0,\infty))\big)$ \ are singular,
 where \ $\cB(C[0,\infty))$ \ denotes the Borel $\sigma$-algebra on \ $C[0,\infty)$.
\ We also note that our technique for the proof of Theorem \ref{THM_singular2} differs from the technique
 of Prakasa Rao \cite[Theorem 5]{Rao_2008}.
Our proof is based on strong consistency of the maximum likelihood estimator of \ $\alpha$,
 \ while the proof of Prakasa Rao \cite[Theorem 5]{Rao_2008} is based on a Baxter type result of Kurchenko
 \cite{Kur} for second order quadratic variations (second order increments) for a fractional Wiener process.
By giving a second proof of  Theorem \ref{THM_singular2}, we also discuss the connections between strong
 consistency of the maximum likelihood estimator of \ $\alpha$, \ Hellinger processes (see, e.g.,
 Jacod and Shiryaev \cite[Chapter IV]{JacShi}) and singularity of induced measures.
Moreover, we study absolute continuity and singularity of probability measures induced by
 processes for which the diffusion coefficients in the SDE \eqref{alpha_W_bridge} are not identically one.
Giving two different proofs, we prove that a so-called dichotomy holds, see Theorem \ref{THM_singular1}.

\section{Preliminaries}\label{Section_Preliminaries}

First we determine the covariance of \ $X_s^{(\alpha)}$ \ and \ $X_t^{(\beta)}$ \ for all
 \ $s,\,t\in[0,T)$ \ and \ $\alpha,$ $\beta\in\RR$.

\begin{Lem}\label{LEMMA_kovariancia}
Let \ $T\in(0,\infty)$ \ be fixed.
For \ $\alpha$, $\beta\in\RR$, \ let us consider the processes \ $(X_t^{(\alpha)})_{t\in[0,T)}$ \
 and \ $(X_t^{(\beta)})_{t\in[0,T)}$ \ given by the SDE \eqref{alpha_W_bridge}.
Then for all \ $s,\,t\in[0,T)$, \ the covariance of \ $X_s^{(\alpha)}$ \ and
 \ $X_t^{(\beta)}$ \ is
 \begin{align}\label{SEGED_COV2}
   \cov(X^{(\alpha)}_s,X^{(\beta)}_t)
      =\begin{cases}
          \frac{(T-s)^\alpha(T-t)^\beta}{1-\alpha-\beta}
            \left(T^{1-\alpha-\beta}-(T-(s\wedge t))^{1-\alpha-\beta}\right)
              & \text{if \ $\alpha+\beta\ne1$, }\\[1mm]
           (T-s)^\alpha(T-t)^\beta
          \ln\left(\frac{T}{T-(s\wedge t)}\right)
               & \text{if \ $\alpha+\beta=1$.}
        \end{cases}
 \end{align}
Especially, for all \ $t\in[0,T)$, \ $X_t^{(\alpha)}$ \ is a normally distributed random variable
 with mean \ $\EE X_t^{(\alpha)}=0$ \ and with variance
 \begin{align*}
    \EE (X_t^{(\alpha)})^2
       =\begin{cases}
           \frac{T}{1-2\alpha}\left(\frac{T-t}{T}\right)^{2\alpha}
               -\frac{T-t}{1-2\alpha}
                   & \text{\quad if \ $\alpha\not=\frac{1}{2}$, }\\
           (T-t)(\ln(T)-\ln(T-t)) & \text{\quad if \ $\alpha=\frac{1}{2}$.}
        \end{cases}
 \end{align*}
\end{Lem}

\noindent{\bf Proof.}
By Bauer \cite[Lemma 48.2]{Bauer} and \eqref{alpha_W_bridge_intrep},
 \ $(X_t^{(\alpha)})_{t\in[0,T)}$ \ is a Gauss process with mean
 \ $\EE X_t^{(\alpha)}=0$, \ $t\in[0,T)$, \ and with variance
 \ $\EE (X_t^{(\alpha)})^2$, \ $t\in[0,T)$.
\ By Proposition 3.2.10 in Karatzas and Shreve \cite{KarShr},
 for all \ $s,\,t\in[0,T)$, \ we have
 \begin{align*}
   \cov(X^{(\alpha)}_s,X^{(\beta)}_t)
    &=\EE(X^{(\alpha)}_sX^{(\beta)}_t)
     =\EE\left(\int_0^s\left(\frac{T-s}{T-u}\right)^\alpha\,\dd B_u
              \int_0^t\left(\frac{T-t}{T-v}\right)^\beta\,\dd B_v\right)\\
    &=(T-s)^\alpha(T-t)^\beta
      \int_0^{s\wedge t}\frac{1}{(T-u)^{\alpha+\beta}}\,\dd u,
 \end{align*}
 and hence we obtain \eqref{SEGED_COV2}.
\proofend

For proving regularity properties of \ $X^{(\alpha)}$, \ we recall a strong law of
 large numbers and a law of the iterated logarithm for continuous local martingales.

Let \ $T\in(0,\infty]$ \ be fixed.
In all what follows, if \ $(M_t)_{t \in [0,T)}$ \ is a continuous local martingale
 satisfying \ $\PP(M_0=0)=1$, \ then \ $(\langle M\rangle_t)_{t\in[0,T)}$ \ denotes
 the quadratic variation of \ $M$.

The following theorem is a modification of Theorem 3.4.6 in Karatzas and Shreve \cite{KarShr}
 (due to Dambis, Dubins and Schwartz), see also Theorem 1.6 in Chapter V in Revuz and
 Yor \cite{RevYor}.
In fact, our next Theorem \ref{DDS} is Exercise 1.18 in Chapter V in Revuz and Yor \cite{RevYor}.

\begin{Thm}\label{DDS}
Let \ $T\in(0,\infty]$ \ be fixed and let \ $\big(\Omega, \cG, (\cG_t)_{t\in[0,T)}, \PP\big)$
 \ be a filtered probability space satisfying the usual conditions.
Let \ $(M_t)_{t \in [0,T)}$ \ be a continuous local martingale with respect to the filtration
 \ $(\cG_t)_{t\in[0,T)}$ \ such that \ $\PP(M_0=0)=1$ \ and
 \ $\PP(\lim_{t \uparrow T} \langle M \rangle_t=\infty)=1$.
\ For each \ $s \in [0,\infty)$, \ define the stopping time
 \[
   \tau_s := \inf \{ t \in [0,T) : \langle M \rangle_t > s \} .
 \]
Then the time-changed process
 \[
   \left( B_s:=M_{\tau_s}, \, \cG_{\tau_s} \right)_{s\geq 0}
 \]
 is a standard Wiener process.
In particular, the filtration \ $(\cG_{\tau_s})_{s\geq 0}$ \ satisfies the usual conditions and
 \[
  \PP\big(M_t=B_{\langle M\rangle_t}\quad\text{for all \ $t\in[0,T)$}\big)=1.
 \]
\end{Thm}

Now we formulate a strong law of large numbers for continuous local martingales.
Compare with L\'epingle \cite[Theoreme 1]{Lep} or with $3^\circ)$ in Exercise 1.16
 in Chapter V in Revuz and Yor \cite{RevYor}.
We note that the above mentioned citations are about continuous local martingales
 with time interval \ $[0,\infty)$, \ but they are also valid for continuous local
 martingales with time interval \ $[0,T)$, \ $T\in(0,\infty)$, \ with appropriate
 modifications in the conditions, see as follows.

\begin{Thm}\label{MSLLN}
Let \ $T\in(0,\infty]$ \ be fixed and let
 \ $\big(\Omega, \cG, (\cG_t)_{t\in[0,T)}, \PP\big)$ \ be a filtered probability space
 satisfying the usual conditions.
Let \ $(M_t)_{t \in [0,T)}$ \ be a continuous local martingale with respect to the
 filtration \ $(\cG_t)_{t\in[0,T)}$ \ such that \ $\PP(M_0=0)=1$ \ and
 \ $\PP(\lim_{t \uparrow T} \langle M \rangle_t = \infty)=1$.
\ Let \ $f : [1, \infty) \to (0, \infty)$ \ be an increasing function such that
 \[
   \int_1^\infty \frac{1}{f(x)^2} \, \dd x < \infty .
 \]
Then
 \[
   \PP\left( \lim_{t \uparrow T} \frac{M_t}{f(\langle M \rangle_t)} = 0 \right)
   = 1 .
 \]
\end{Thm}

\ Now we present a law of the iterated logarithm for continuous local martingales.
Compare with Exercise 1.15 in Chapter V in Revuz and Yor \cite{RevYor}.
We note that the above mentioned citation is about continuous local martingales with time
 interval \ $[0,\infty)$, \ but the result is also valid for continuous local martingales with time
 interval \ $[0,T)$, \ $T\in(0,\infty)$, \ with appropriate modifications in the conditions,
 see as follows.

\begin{Thm}\label{MLIL}
Let \ $T\in(0,\infty]$ \ be fixed and let \ $\big(\Omega, \cG, (\cG_t)_{t\in[0,T)}, \PP\big)$
 \ be a filtered probability space satisfying the usual conditions.
Let \ $(M_t)_{t \in [0,T)}$ \ be a continuous local martingale with respect to the filtration
 \ $(\cG_t)_{t\in[0,T)}$ \ such that \ $\PP(M_0=0)=1$ \ and
 \ $\PP(\lim_{t \uparrow T} \langle M \rangle_t = \infty)=1$.
\ Then
 \[
   \PP\left( \limsup_{t \uparrow T}
              \frac{M_t}{\sqrt{2 \langle M \rangle_t \ln(\ln \langle M \rangle_t)}}
        = 1 \right)
   =
    \PP\left( \liminf_{t \uparrow T}
               \frac{M_t}{\sqrt{2 \langle M \rangle_t \ln(\ln \langle M \rangle_t)}}
        = - 1 \right)
   = 1 .
 \]
\end{Thm}

Theorem \ref{MLIL} simply follows from Theorem \ref{DDS} and the law of the
 iterated logarithm for a standard Wiener process (see, e.g.,
 Karatzas and Shreve \cite[Theorem 2.9.23]{KarShr}).

\section{Sample paths properties}\label{Section_sample_paths}

The following Lemma \ref{LEMMA5} can be considered as a generalization of Lemma 5.6.9 in
 Karatzas and Shreve \cite{KarShr}.
Namely, our result with \ $\alpha=1$ \ gives back Lemma 5.6.9 in Karatzas and Shreve \cite{KarShr}.
Furthermore, Lemma \ref{LEMMA5} can be also considered as a generalization of Corollary 4.4
 in Becker-Kern \cite{Bec} in the $1$-dimensional Brownian case.
Namely, by \eqref{alpha_W_bridge_intrep}, \ $(X_t^{(\alpha)})_{t\in[0,T)}$ \ coincides
 with the process \ $(U_t)_{t\in[0,T)}$ \ defined in (4.7) in Becker-Kern \cite{Bec} in the
 $1$-dimensional Brownian case for all \ $\alpha>0$.
\ In this case Becker-Kern proved that \ $U_t$ \ converges in probability to
 \ $0$ \ as \ $t\uparrow T$, \ while we prove convergence with probability one.
For historical fidelity, we remark that something similar to the statement of our
 Lemma \ref{LEMMA5} is stated on page 1023 in Mansuy \cite{Man} but without any proof.

\begin{Lem}\label{LEMMA5}
Let \ $T\in(0,\infty)$ \ and \ $\alpha>0$ \ be fixed, and let \ $(B_t)_{t\geq 0}$ \ be
 a $1$-dimensional standard Wiener process.
The process \ $(Y_t^{(\alpha)})_{t\in[0,T]}$ \ defined by
 \[
  Y_t^{(\alpha)}:=\begin{cases}
            \int_0^t\left(\frac{T-t}{T-s}\right)^\alpha\,\dd B_s
                & \text{\quad if \ $t\in[0,T)$,}\\
                 0 & \text{\quad if \ $t=T$,}
        \end{cases}
 \]
 is a centered Gauss process with almost surely continuous paths.
\end{Lem}

\noindent{\bf Proof.}
By Bauer \cite[Lemma 48.2]{Bauer}, \ $(Y_t^{(\alpha)})_{t\in[0,T]}$ \
 is a centered Gauss process.
To prove almost surely continuity, we follow the method of the proof of
 Lemma 5.6.9 in Karatzas and Shreve \cite{KarShr}.
For all \ $t\in[0,T)$ \ and \ $\alpha\in\RR$, \ let
 \[
   M^{(\alpha)}_t:=\int_0^t\frac{1}{(T-s)^\alpha}\,\dd B_s.
 \]
Then \ $(M_t^{(\alpha)})_{t\in[0,T)}$ \ is a continuous, square-integrable martingale
 with respect to the filtration induced by \ $B$ \ and with quadratic variation
 \begin{align}\label{SEGED_quadratic}
   \langle M^{(\alpha)}\rangle_t
      :=\int_0^t\frac{1}{(T-s)^{2\alpha}}\,\dd s
       =\begin{cases}
         \frac{T^{1-2\alpha}}{1-2\alpha}\left(1-\left(1-\frac{t}{T}\right)^{1-2\alpha}\right)
             & \text{\quad if \ $\alpha\not=\frac{1}{2}$, }\\
         -\ln\left(1-\frac{t}{T}\right) & \text{\quad if \ $\alpha=\frac{1}{2}$, }
        \end{cases}
         \quad t\in[0,T).
 \end{align}
Then
 \begin{align}\label{SEGED11}
   \lim_{t\uparrow T}\langle M^{(\alpha)}\rangle_t
      =\begin{cases}
          \infty & \text{\quad if \ $\alpha\geq\frac{1}{2}$, }\\
          \frac{T^{1-2\alpha}}{1-2\alpha}
                  & \text{\quad if \ $\alpha<\frac{1}{2}$.}
       \end{cases}
 \end{align}
Hence in case of \ $\alpha\geq \frac{1}{2}$, \ Theorem \ref{DDS}
 implies that there exists a standard $1$-dimensional Wiener process \ $(W_t)_{t\geq 0}$
 \ on \ $(\Omega,\cA,\PP)$ \ such that
 \[
  \PP\big(M^{(\alpha)}_t=W_{\langle M^{(\alpha)}\rangle_t}
      \quad\text{for all \ $t\in[0,T)$}\big)=1.
 \]
First we consider the case of \ $\alpha>\frac{1}{2}$.
\ Let us define the function \ $f_\alpha:[1,\infty)\to(0,\infty)$ \ by
 \ $f_\alpha(x):=x^{\alpha/(2\alpha-1)}$, \ $x\geq 1$.
\ Then \ $f_\alpha$ \ is strictly monotone increasing and
 \[
   \int_1^\infty \frac{1}{f_\alpha(x)^2} \, \dd x
   = \int_1^\infty x^{ - 2 \alpha / (2 \alpha - 1)} \, \dd x
   = 2 \alpha - 1
   < \infty ,
 \]
 hence we may apply Theorem \ref{MSLLN} and then we obtain
 \[
   \PP\left( \lim_{t \uparrow T}
               \frac{M^{(\alpha)}_t}{f_\alpha(\langle M^{(\alpha)} \rangle_t)} = 0
       \right)
   = 1 , \qquad \alpha>\frac{1}{2}.
 \]
We have
 \[
   Y_t^{(\alpha)}
   = (T-t)^\alpha M_t^{(\alpha)}
   = (T-t)^\alpha f_\alpha(\langle M^{(\alpha)} \rangle_t) \,
     \frac{M_t^{(\alpha)}}{f_\alpha(\langle M^{(\alpha)} \rangle_t)} ,
 \]
 where \ $t\in[0,T)$ \ is such that \ $\langle M^{(\alpha)} \rangle_t\geq 1$.
\ Here
 \[
   (T-t)^\alpha f_\alpha(\langle M^{(\alpha)} \rangle_t)
   \leq (T-t)^\alpha
        f_\alpha\left( \frac{1}{2\alpha-1} \, \frac{1}{(T-t)^{2\alpha-1}}\right)
   = (2 \alpha - 1)^{- \alpha / (2 \alpha - 1)},
 \]
 where \ $t \in [0,T)$ \ is such that \ $\langle M^{(\alpha)} \rangle_t\geq 1$.
\ Hence we conclude
 \ $\PP\left( \lim\limits_{t \uparrow T} Y^{(\alpha)}_t = 0 \right) = 1$.

\noindent Now we consider the case of \ $\alpha = \frac{1}{2}$.
\ Let us define the function \ $f_{1/2}(x) := \ee^{x / 2}$, \ $x\in [1, \infty)$.
\ Then \ $f_{1/2}$ \ is strictly monotone increasing and
 \[
   \int_1^\infty \frac{1}{f_{1/2}(x)^2} \, \dd x
   = \int_1^\infty  \ee^{-x} \, \dd x
   = \ee^{-1}
   < \infty ,
 \]
 hence we may apply Theorem \ref{MSLLN} and then we obtain
 \[
   \PP\left( \lim_{t \uparrow T}
               \frac{M^{(1/2)}_t}{f_{1/2}(\langle M^{(1/2)} \rangle_t)} = 0
       \right)
   = 1 .
 \]
We have
 \[
   Y_t^{(1/2)}
   = (T-t)^{1/2} M_t^{(1/2)}
   = (T-t)^{1/2} f_{1/2}(\langle M^{(1/2)} \rangle_t) \,
     \frac{M_t^{(1/2)}}{f_{1/2}(\langle M^{(1/2)} \rangle_t)} ,
 \]
 where \ $t \in [0,T)$ \ is such that \ $\langle M^{(1/2)} \rangle_t\geq 1$.
\ Here
 \[
   (T-t)^{1/2} f_{1/2}(\langle M^{(1/2)} \rangle_t)
    = (T-t)^{1/2} \exp\left\{ \frac{1}{2} \, \ln\left(\frac{T}{T-t}\right)\right\}
    = T^{1/2},
 \]
 where \ $t \in [0,T)$ \ is such that \ $\langle M^{(1/2)} \rangle_t\geq 1$.
\ Hence we conclude
 \ $\PP\left( \lim\limits_{t \uparrow T} Y^{(1/2)}_t = 0 \right) = 1$.

\noindent Finally, we consider the case of \ $0<\alpha<\frac{1}{2}$.
\ Using \eqref{SEGED11} we have Proposition 1.26 in Chapter IV and
 Proposition 1.8 in Chapter V in Revuz and Yor \cite{RevYor} imply that
 the limit \ $M^{(\alpha)}_T:=\lim_{t\uparrow T}M^{(\alpha)}_t$ \ exists
 almost surely.
Since
 \[
  Y_t^{(\alpha)}
        =(T-t)^\alpha M^{(\alpha)}_t,
        \qquad t\in[0,T),
 \]
 we get \ $\lim_{t\uparrow T}Y_t^{(\alpha)}=0$ \ almost surely.
\proofend

By Lemma \ref{LEMMA5}, we can say that in case of \ $\alpha>0$, \ the process
 \ $X^{(\alpha)}$ \ has an almost surely continuous extension.
In the later Remark \ref{REMARK2} we examine the possibility of such an almost surely continuous
 extension of \ $(X_t^{(\alpha)})_{t\in[0,T)}$ \ in case of \ $\alpha\leq 0$.

Now we prove some results about the asymptotic behavior of \ $X^{(\alpha)}_t$
 \ as \ $t \uparrow T$.
\ Theorem \ref{MLIL} has the following consequences on \ $X^{(\alpha)}$.

\begin{Thm}\label{THM_path1}
If \ $\alpha>\frac{1}{2}$, \ then
  \begin{align}\label{SEGED_regular2}
   \PP\left( \limsup_{t \uparrow T}
               \frac{X_t^{(\alpha)}}
                    {\sqrt{\frac{2(T-t)}{2\alpha-1}\ln\left(\ln \frac{1}{T-t}\right)}} = 1
        \right)
   =
   \PP\left(\liminf_{t \uparrow T}
               \frac{X_t^{(\alpha)}}
                    {\sqrt{\frac{2(T-t)}{2\alpha-1}\ln\left(\ln \frac{1}{T-t}\right)}} = - 1 \right)
   = 1.
 \end{align}
Especially,
 \begin{align}\label{SEGED_regular1}
   \PP\left(\limsup_{t\uparrow T}\frac{X_t^{(\alpha)}}{(T-t)^\alpha}=\infty\right)
    =\PP\left(\liminf_{t\uparrow T}\frac{X_t^{(\alpha)}}{(T-t)^\alpha}=-\infty\right)
    =1.
 \end{align}
\end{Thm}

\noindent{\bf Proof.}
With the notation introduced in the proof of Lemma \ref{LEMMA5}, we have
 \ $\frac{X_t^{(\alpha)}}{(T-t)^\alpha}=M_t^{(\alpha)}$, \ $t\in[0,T)$,
 \ and the quadratic variation \ $\langle M^{(\alpha)}\rangle_t$, \ $t\in[0,T)$, \
 of the continuous martingale \ $\frac{X_t^{(\alpha)}}{(T-t)^\alpha}$, \ $t\in[0,T)$,
 \ is given in \eqref{SEGED_quadratic}.
Using \eqref{SEGED11} we have \ $\lim_{t\uparrow T}\langle M^{(\alpha)}\rangle_t=\infty$.
\ Then, by Theorem \ref{MLIL}, in case of \ $\alpha>\frac{1}{2}$ \ we get
 \begin{align*}
   \PP\left(\limsup_{t\uparrow T}
          \frac{X_t^{(\alpha)}}
           {\sqrt{D_{\alpha,T}(t)}}
           =1
       \right)
   =\PP\left(\liminf_{t\uparrow T}
          \frac{X_t^{(\alpha)}}
           {\sqrt{D_{\alpha,T}(t)}}
             =-1
        \right)
   =1,
 \end{align*}
 where for all \ $t\in[0,T)$,
 \[
   D_{\alpha,T}(t)
    :=2(T-t)^{2\alpha}\frac{T^{1-2\alpha}}{1-2\alpha}
         \left(1-\left(1-\frac{t}{T}\right)^{1-2\alpha}\right)
         \ln\!\left(\!\ln\!\left(\!\frac{T^{1-2\alpha}}{1-2\alpha}
         \left(1-\left(1-\frac{t}{T}\right)^{1-2\alpha}\!\right)\!\right)\!\right).
 \]
Hence to prove \eqref{SEGED_regular2} it is enough to check that
 \begin{align*}
   \lim_{t\uparrow T}
     \frac{2(T-t)^{2\alpha}\frac{T^{1-2\alpha}}{1-2\alpha}
               \left(1-\left(1-\frac{t}{T}\right)^{1-2\alpha}\right)
                \ln\left(\ln\left(\frac{T^{1-2\alpha}}{1-2\alpha}
               \left(1-\left(1-\frac{t}{T}\right)^{1-2\alpha}\right)\right)\right)}
           {\frac{2(T-t)}{2\alpha-1}\ln\left(\ln \frac{1}{T-t}\right)}
       =1.
 \end{align*}
This is satisfied, since, by using L'Hospital's rule twice, we get
 \begin{align*}
   &\lim_{t\uparrow T}
       \frac{\ln\left(\ln\left(\frac{T^{1-2\alpha}}{1-2\alpha}
               \left(1-\left(1-\frac{t}{T}\right)^{1-2\alpha}\right)\right)\right)}
              {\ln\left(\ln \frac{1}{T-t}\right)}
       =\lim_{t\uparrow T}
          \frac{(T-t)^{1-2\alpha}-T^{1-2\alpha}}{(T-t)^{1-2\alpha}}
       =1.
 \end{align*}
Using \eqref{SEGED_regular2} and the decomposition
 \[
   \frac{X_t^{(\alpha)}}{(T-t)^\alpha}
     =\frac{X_t^{(\alpha)}}
        {\sqrt{\frac{2(T-t)}{2\alpha-1}\ln\left(\ln \frac{1}{T-t}\right)}}
        \sqrt{\frac{2(T-t)^{1-2\alpha}}{2\alpha-1}\ln\left(\ln \frac{1}{T-t}\right)},
      \quad t\in[0,T),
 \]
 we have \eqref{SEGED_regular1}.
\proofend

The next theorem is about the limit behavior of \ $X_t^{(1/2)}$ \ as
\ $t\uparrow T$. \

\begin{Thm}\label{THM_path2}
We have
 \begin{align}\label{SEGED_regular4}
  \begin{split}
    &\PP\left(  \limsup_{t \uparrow T}
               \frac{X_t^{(1/2)}}
                    {\sqrt{2 (T-t) \left(\ln \frac{1}{T-t} \right)
                           \left(\ln \ln \ln \frac{1}{T-t} \right)}} = 1
      \right)\\[2mm]
   &\qquad
    =\PP\left(   \liminf_{t \uparrow T}
               \frac{X_t^{(1/2)}}
                    {\sqrt{2 (T-t) \left(\ln \frac{1}{T-t} \right)
                           \left(\ln \ln \ln \frac{1}{T-t} \right) }} = - 1
           \right)
   = 1.
  \end{split}
 \end{align}
Especially,
\begin{align*}
 \PP\left(\limsup_{t\uparrow T}\frac{X_t^{(1/2)}}{\sqrt{T-t}}=\infty\right)
    =\PP\left(\liminf_{t\uparrow T}\frac{X_t^{(1/2)}}{\sqrt{T-t}}=-\infty\right)
    =1.
\end{align*}
\end{Thm}

\noindent{\bf Proof.}
With the notation introduced in the proof of Lemma \ref{LEMMA5}, we have
 \ $\frac{X_t^{(1/2)}}{\sqrt{T-t}}=M_t^{(1/2)}$,  $t\in[0,T)$,  \
 and the quadratic variation \ $\langle M^{(1/2)}\rangle_t$, \ $t\in[0,T)$, \ of
 the continuous martingale \ $\frac{X_t^{(1/2)}}{\sqrt{T-t}}$, \ $t\in[0,T)$,
 \ is given in \eqref{SEGED_quadratic}.
Using \eqref{SEGED11} we have \ $\lim_{t\uparrow T}\langle M^{(1/2)}\rangle_t=\infty$.
\ Then, by Theorem \ref{MLIL}, we get
 \begin{align*}
  &\PP\left(  \limsup_{t \uparrow T}
               \frac{X_t^{(1/2)}}
                    {\sqrt{2 (T-t) \left(\ln \frac{T}{T-t} \right)
                           \left(\ln \ln \ln \frac{T}{T-t} \right)}} = 1
      \right)\\[2mm]
   &\qquad
    =\PP\left(   \liminf_{t \uparrow T}
               \frac{X_t^{(1/2)}}
                    {\sqrt{2 (T-t) \left(\ln \frac{T}{T-t} \right)
                           \left(\ln \ln \ln \frac{T}{T-t} \right) }} = - 1
           \right)
   = 1,
 \end{align*}
 which yields \eqref{SEGED_regular4}.
\proofend

The next theorem is about the limit behavior of \ $X_t^{(\alpha)}$ \ as \ $t\uparrow T$ \
in case of \ $\alpha<\frac{1}{2}$. \

\begin{Thm}\label{THM_path3}
If \ $\alpha < \frac12$,  \ then
\begin{align}\label{SEGED_regular5_spec}
   &\PP\left( \lim_{t \uparrow T} \frac{X_t^{(\alpha)}} {(T-t)^\alpha} = M_T^{(\alpha)}\right)=1,
 \end{align}
 where \ $M_T^{(\alpha)}$ \ is a normally distributed random variable with mean 0
 and with variance \ $\frac{T^{1-2\alpha}}{1-2\alpha}$.
\ Consequently,
 \begin{align}
   \label{SEGED_regular5}
   &\PP\left( \lim_{t \uparrow T} \frac{X_t^{(\alpha)}} {(T-t)^\beta} = 0 \right)
     = 1 \qquad \text{for all \ $\beta < \alpha$, }\\\label{SEGED_regular6}
   &\PP\left( \lim_{t \uparrow T} \frac{X_t^{(\alpha)}} {(T-t)^\beta} = -\infty
        \right)
     =
    \PP\left( \lim_{t \uparrow T} \frac{X_t^{(\alpha)}} {(T-t)^\beta} = \infty
        \right)
      = \frac{1}{2}
      \qquad \text{for all \ $\beta > \alpha$.}
 \end{align}
\end{Thm}

\noindent{\bf Proof.}
By \eqref{alpha_W_bridge_intrep}, using the notations introduced in the proof of Lemma \ref{LEMMA5},
 we get
 \[
  M^{(\alpha)}_t
   =\frac{X^{(\alpha)}_t}{(T-t)^\alpha}
   =\int_0^t\frac{1}{(T-s)^\alpha}\,\dd B_s,
   \qquad t\in[0,T).
 \]
By \eqref{SEGED11}, since \ $\alpha<\frac{1}{2}$,
 \[
  \lim_{t\uparrow T}\langle M^{(\alpha)}\rangle_t
     =\frac{T^{1-2\alpha}}{1-2\alpha}<\infty,
 \]
 and hence Proposition 1.26 in Chapter IV and Proposition 1.8 in Chapter V
 in Revuz and Yor \cite{RevYor} imply that the limit
 \ $M^{(\alpha)}_T:=\lim_{ t\uparrow T}M^{(\alpha)}_t$ \ exists
 almost surely.
Using that \ $M^{(\alpha)}_t$ \ is normally distributed with mean \ $0$
 \ and with variance \ $\langle M^{(\alpha)}\rangle_t$ \ for all \ $t\in[0,T)$,
 \ we have the random variable \ $M^{(\alpha)}_T$ \ is also normally distributed with
 mean \ $0$ \ and with variance \ $\frac{T^{1-2\alpha}}{1-2\alpha}$.
\ Indeed, normally distributed random variables can converge in distribution only to a
 normally distributed random variable, by continuity theorem, see, e.g., page 304 in Shiryaev
 \cite{Shi}.
This implies \eqref{SEGED_regular5_spec}.
Hence for all \ $\alpha, \beta\in\RR$, \ we get
 \begin{align}\label{SEGED_regular7}
   \frac{X^{(\alpha)}_t}{(T-t)^\beta}
     =(T-t)^{\alpha-\beta}M^{(\alpha)}_t,
     \qquad t\in[0,T).
 \end{align}
If \ $\beta<\alpha$, \ then using \eqref{SEGED_regular7} and that
 \ $\PP(\lim_{t\uparrow T} M^{(\alpha)}_t=M^{(\alpha)}_T)=1$, \ we get \eqref{SEGED_regular5}.
If \ $\beta>\alpha$, \ using that \ $\PP(M^{(\alpha)}_T=0)=0$, \ we have \eqref{SEGED_regular7}
 implies that
 \[
   \PP\left(\lim_{t\uparrow T}\frac{X^{(\alpha)}_t}{(T-t)^\beta}
           \in\{-\infty,\infty\}\right)=1.
 \]
Since \ $\PP(M^{(\alpha)}_T>0)=\PP(M^{(\alpha)}_T<0)=\frac{1}{2}$, \ we get \eqref{SEGED_regular6}.
\proofend

\begin{Rem}\label{REMARK2}
In case of \ $\alpha=0$, \ the process \ $(X_t^{(0)})_{t\in[0,T)}$ \ is a standard
 Wiener process and hence it can be extended to an almost surely continuous process
 \ $(Y_t^{(0)})_{t\in[0,T]}$ \ with the definition \ $Y_t^{(0)}:=B_T$.
\ In case of \ $\alpha<0$, \ there does not exist an almost surely continuous process
 \ $(Y_t^{(\alpha)})_{t\in[0,T]}$ \ such that \ $\PP(X_t^{(\alpha)}=Y_t^{(\alpha)})=1$
 \ for all \ $t\in[0,T)$.
\ Indeed, by \eqref{SEGED_regular6}, we get
 \[
   \PP\left( \lim_{t \uparrow T} X_t^{(\alpha)} =  -\infty
        \right)
     =
    \PP\left( \lim_{t \uparrow T} X_t^{(\alpha)} =  \infty
        \right)
      = \frac{1}{2}
      \qquad \text{for all \ $\alpha <0$.}
 \]
\proofend
\end{Rem}

\section{Singularity of induced measures}\label{Section_measures}

For probability measures \ $\PP_1$ \ and \ $\PP_2$ \ on a measurable space \ $(\Omega,\cG)$,
 \ equivalence and singularity of them will be denoted by \ $\PP_1\sim \PP_2$ \ and \ $\PP_1\perp \PP_2$,
 \ respectively.

Using that for all \ $\alpha\in\RR$, \ the process \ $X^{(\alpha)}$
 \ has continuous paths (by the definition of strong solution,
 see, e.g., Jacod and Shiryaev \cite[Definition 2.24, Chapter III]{JacShi}),
 we have
 \begin{align}\label{SEGED13}
   \PP\left(\int_0^t\frac{(X_u^{(\alpha)})^2}{(T-u)^2}\,\dd u<\infty\right)=1,
       \qquad \forall\;\alpha\in\RR,\;\;\; t\in[0,T).
 \end{align}
For all \ $\alpha\in\RR$ \ and \ $t\in(0,T)$, \ let \ $\PP_{X^{(\alpha)},\,t}$
 \ denote the law of the process \ $(X_s^{(\alpha)})_{s\in[0,t]}$ \ on
 \ $\big(C[0,t],\cB(C[0,t])\big)$, \ where \ $\cB(C[0,t])$ \ denotes the
 Borel $\sigma$-algebra on \ $C[0,t]$.
 \ Using Theorem 7.20 in Liptser and Shiryaev \cite{LipShiI} and \eqref{SEGED13},
 we get \ $\PP_{X^{(\alpha)},\,t}\sim \PP_{X^{(0)},\,t}$ \ and
 \begin{align}\label{SEGED16}
   \frac{\dd\,\PP_{X^{(\alpha)},\,t}}{\dd\,\PP_{X^{(0)},\,t}}(X^{(\alpha)}\vert_{[0,t]})
       &=\exp\left\{-\alpha\int_0^t\frac{X_u^{(\alpha)}}{T-u}\,\dd X_u^{(\alpha)}
          -\frac{\alpha^2}{2}\int_0^t\frac{(X_u^{(\alpha)})^2}{(T-u)^2}\,\dd u\right\}.
 \end{align}
Here \ $\PP_{X^{(0)},\,t}$ \ is nothing else but the Wiener measure on \ $\big(C[0,t],\cB(C[0,t])\big)$.

We recall that for all \ $t\in(0,T)$, \ the maximum likelihood estimator (MLE)
 \ $\halpha_t^{(X^{(\alpha)})}$ \ of the parameter \ $\alpha$ \ based on the
 observation \ $(X_s^{(\alpha)})_{s \in [0,\,t]}$ \ is defined by
 \[
   \halpha_t^{(X^{(\alpha)})}
    :=\argmax_{\alpha\in\RR}
       \ln\left(\frac{\dd \PP_{X^{(\alpha)},\,t}}{\dd \PP_{X^{(0)},\,t}}
                 \left(X^{(\alpha)}\big\vert_{[0,t]}\right)\right).
 \]
By \eqref{SEGED13} and \eqref{SEGED16}, for all \ $t\in(0,T)$, \ there exists a unique
 MLE \ $\halpha_t^{(X^{(\alpha)})}$ \ of the parameter \ $\alpha$ \ based on the observation
 \ $(X_s^{(\alpha)})_{s \in [0,\,t]}$ \ given by
 \[
   \halpha_t^{(X^{(\alpha)})}
   = -\frac{\int_0^t
            \frac{X_s^{(\alpha)}}{T-s} \, \dd X_s^{(\alpha)}}
          {\int_0^t
            \frac{(X_s^{(\alpha)})^2}{(T-s)^2} \, \dd s},
   \qquad  t\in(0,T).
 \]
To be more precise, by \eqref{SEGED13}, for all \ $t\in(0,T)$, \ the MLE \ $\halpha_t^{(X^{(\alpha)})}$
 \ exists $\PP$-almost surely. As a special case of Theorem 3.12 in Barczy and Pap \cite{BarPap2},
 the MLE of \ $\alpha$ \ is strongly consistent, i.e.,
 \begin{align}\label{STRONG_CONSISTENCY}
   \PP\big(\lim_{t\uparrow T}\halpha_t^{(X^{(\alpha)})}
                    =\alpha\big)=1,\qquad \alpha\in\RR.
 \end{align}

For all \ $\alpha\in\RR$, \ let \ $\PP_{X^{(\alpha)},\,T}^T$ \ be the law of the process
 \ $(X_t^{(\alpha)})_{t\in[0,T)}$ \ given by the SDE \eqref{alpha_W_bridge}
 on \ $(C[0,T),\cB(C[0,T)))$.

\begin{Thm}\label{THM_singular2}
For all \ $\alpha, \beta\in\RR$, \ $\alpha\ne\beta$, \ we have
 \ $\PP_{X^{(\alpha)},\,T}^T \perp \PP_{X^{(\beta)},\,T}^T$.
\ In other words, the laws of the processes \ $(X_t^{(\alpha)})_{t\in[0,T)}$ \ and
 \ $(X_t^{(\beta)})_{t\in[0,T)}$ \ on \ $(C[0,T),\cB(C[0,T)))$ \ are singular for
 all \ $\alpha,\,\beta\in\RR$, \ $\alpha\ne\beta$.
\end{Thm}

\noindent{\bf Proof.}
First we check that for all \ $\alpha\in\RR$ \ and \ $t\in[0,T)$,
 \begin{align}\label{SEGED_LAPLACE1}
     \int_0^t\frac{X^{(\alpha)}_s}{T-s} \,\dd X_s^{(\alpha)}
        =\frac{1}{2}\left(\frac{(X^{(\alpha)}_t)^2}{T-t}
          -\int_0^t\frac{(X^{(\alpha)}_s)^2}{(T-s)^2}\,\dd s
                 -\ln\left(\frac{T}{T-t}\right)\right).
 \end{align}
By It\^{o}'s rule (see, e.g., Liptser and Shiryaev \cite[Theorem 4.4]{LipShiI}), we get
 \begin{align}\label{SEGED_ITO}
  \begin{split}
  \dd\left(\frac{X^{(\alpha)}_t}{T-t}\right)
     &=\frac{X^{(\alpha)}_t}{(T-t)^2}\,\dd t
        +\frac{1}{T-t}\,\dd X^{(\alpha)}_t\\
     &=\frac{X^{(\alpha)}_t}{(T-t)^2}\,\dd t
       -\alpha\frac{X^{(\alpha)}_t}{(T-t)^2}\,\dd t
       +\frac{1}{T-t}\,\dd B_t,
       \quad t\in[0,T).
   \end{split}
 \end{align}
Now we verify that \ $(X_t^{(\alpha)})_{t\in[0,T)}$ \ and
 \ $\left(\frac{X^{(\alpha)}_t}{T-t}\right)_{t\in[0,T)}$
 \ are continuous semimartingales adapted to the filtration induced by \ $B$.
\ Consider the decomposition
 \[
  X_t^{(\alpha)}
     =(T-t)^\alpha
      \int_0^t\frac{1}{(T-s)^\alpha}\,\dd B_s,
      \qquad t\in[0,T).
 \]
Here the deterministic function \ $(T-t)^\alpha$, \ $t\in[0,T)$,
 \ is monotone and hence has a finite variation over each finite interval of
 \ $[0,T)$, \ and then, by Jacod and Shiryaev \cite[Proposition 4.28, Chapter I]{JacShi},
 it is a semimartingale. Since
 \[
    \int_0^t\frac{1}{(T-s)^\alpha}\,\dd B_s,
      \qquad t\in[0,T),
 \]
 is a martingale with respect to the filtration induced by \ $B$,
 \ using Theorem 4.57 in Chapter I in Jacod and Shiryaev \cite{JacShi} with the function
 \ $f(x,y):=xy$, \ $x,y\in\RR$, \ we have \ $(X_t^{(\alpha)})_{t\in[0,T)}$
 \ is a continuous semimartingale adapted to the filtration induced by \ $B$.
\ Similarly as above, using that \ $\frac{1}{T-t}$, \ $t\in[0,T)$,
 \ is continously differentiable, and hence has a finite variation over each finite interval of
 \ $[0,T)$, \ one can get
 \ $\left(\frac{X^{(\alpha)}_t}{T-t}\right)_{t\in[0,T)}$
 \ is a continuous semimartingale adapted to the filtration induced by \ $B$.
\ Moreover, by \eqref{SEGED_ITO}, the cross-variation process of the continuous martingale
 parts of the processes \ $(X^{(\alpha)}_t)_{t\in[0,T)}$ \ and
 \ $\left(\frac{X^{(\alpha)}_t}{T-t}\right)_{t\in[0,T)}$ \ equals
 \[
        \int_0^t \frac{1}{T-s}\,\dd s
           =\ln\left(\frac{T}{T-t}\right), \qquad t\in[0,T).
 \]
Hence, by integration by parts formula (see, e.g., Karatzas and Shreve \cite[page 155]{KarShr}),
 we have for all \ $t\in[0,T)$,
 \begin{align*}
  \int_0^t\frac{X^{(\alpha)}_s}{T-s} \,\dd X_s^{(\alpha)}
     &=\frac{X^{(\alpha)}_t}{T-t}X_t^{(\alpha)}
       -\int_0^tX^{(\alpha)}_s\,\dd \left(\frac{X^{(\alpha)}_s}{T-s}\right)
       -\ln\left(\frac{T}{T-t}\right)\\
     &=\frac{(X^{(\alpha)}_t)^2}{T-t}
       -\int_0^t\frac{(X^{(\alpha)}_s)^2}{(T-s)^2} \,\dd s
       -\int_0^t\frac{X^{(\alpha)}_s}{T-s}
              \,\dd X^{(\alpha)}_s
       -\ln\left(\frac{T}{T-t}\right),
 \end{align*}
 which yields \eqref{SEGED_LAPLACE1}.
Hence \ $\halpha_t^{(X^{(\alpha)})}=A_t(X^{(\alpha)})$ \ for all \ $t\in(0,T)$, \ where
 \ $A_t:C[0,T)\setminus\{0\}\to\RR$, \ defined by
 \[
  A_t(x)
    :=\frac{-\frac{x(t)^2}{T-t}
             +\int_0^t\frac{x(s)^2}{(T-s)^2}\,\dd s
             +\ln\left(\frac{T}{T-t}\right)}
        {2\int_0^t \frac{x(s)^2}{(T-s)^2} \, \dd s},
     \qquad x\in C[0,T)\setminus\{0\},\quad t\in(0,T).
 \]
For all \ $\alpha\in\RR$, \ let us introduce the following subset of \ $C[0,T)$,
 \begin{align*}
    S_\alpha:=\left\{
               x\in C[0,T)\setminus\{0\} :
               \lim_{t\uparrow T}A_t(x)
                =\alpha
                \right\}.
 \end{align*}
We check that \ $S_\alpha\in\cB(C[0,T))$.
\ By Problem 2.4.1 in Karatzas and Shreve \cite{KarShr}, under the metric
 \[
  \rho(x,y):=\sum_{n=1}^\infty\frac{1}{2^n}\sup_{u\in[0,n]}(\vert x(\Psi(u))-y(\Psi(u))\vert\wedge 1),
   \qquad x,y\in C[0,T),
 \]
 the set \ $ C[0,T)$ \ is a complete, separable metric space, where
 \ $\Psi:[0,\infty)\mapsto[0,T)$, \ $\Psi(u):=\frac{2T}{\pi}\arctan(u)$, \ $u\geq0$.
\ For all \ $t\in[0,T)$, \ let \ $L_t:C[0,T)\to\RR$, \
 \[
   L_t(x):=\int_0^t\frac{x(s)^2}{(T-s)^2}\,\dd s,
       \qquad x\in C[0,T).
 \]
Let \ $x\in C[0,T)$ \ be fixed.
We show that for all \ $t\in[0,T)$, \ $L_t$ \ is continuous at the point \ $x\in C[0,T)$.
\ Indeed, for all \ $y\in C[0,T)$, \ we have
 \begin{align*}
  \vert L_t(x)-L_t(y)\vert
     &\leq \sup_{s\in[0,t]}\Big(\vert x(s)+y(s)\vert \vert x(s)-y(s)\vert\Big)
          \int_0^t \frac{1}{(T-s)^2}\,\dd s\\
     &\leq \sup_{s\in[0,t]}\Big((2\vert x(s)\vert+\vert y(s)- x(s)\vert) \vert x(s)-y(s)\vert\Big)
          \int_0^t \frac{1}{(T-s)^2}\,\dd s.
 \end{align*}
If \ $y\in C[0,T)$ \ is such that \ $\delta:=\sup_{s\in[0,t]}\vert y(s)- x(s)\vert<1$ \ and
 \ $n_0\in\NN$ \ is such that \ $n_0>\Psi^{-1}(t)$, \ then
 \[
   \rho(x,y)\geq \frac{1}{2^{n_0}}\sup_{u\in[0,n_0]}(\vert y(\Psi(u))-x(\Psi(u))\vert\wedge 1)
            \geq \frac{1}{2^{n_0}}\sup_{u\in[0,\Psi^{-1}(t)]}\vert y(\Psi(u))-x(\Psi(u))\vert
             = \frac{\delta}{2^{n_0}},
 \]
 and hence
 \begin{align*}
   \vert L_t(x)-L_t(y)\vert
      \leq \delta \big(1+2\sup_{s\in[0,t]}\vert x(s)\vert \big) \int_0^t \frac{1}{(T-s)^2}\,\dd s
      \leq K(t)\rho(x,y),
 \end{align*}
 where \ $K(t):=2^{n_0}\big(1+2\sup_{s\in[0,t]}\vert x(s)\vert \big) \int_0^t \frac{1}{(T-s)^2}\,\dd s$,
 \ which yields the continuity of \ $L_t$ \ at \ $x$.
\ Consequently, \ $A_t$ \ is continuous for all \ $t\in(0,T)$.
\ Consider the decomposition
 \begin{align*}
  S_\alpha
    &=\bigcap_{\varepsilon>0}\bigcup_{t\in[0,T)}\bigcap_{s\in[t,T)}
      \Big\{ x\in C[0,T)\setminus\{0\} : \vert A_s(x)-\alpha\vert\leq \varepsilon\Big\}\\
    &=\bigcap_{n=1}^\infty\bigcup_{m=1}^\infty \bigcap_{s\in[T-\frac{1}{m},T)\cap\QQ_+}
      \left\{ x\in C[0,T)\setminus\{0\} : \vert A_s(x)-\alpha\vert\leq \frac{1}{n}\right\} ,
 \end{align*}
 where \ $\QQ_+$ \ denotes the set of positive rational numbers.
Since \ $A_s$ \ is continuous for all \ $s\in(0,T)$, \ we have
 \[
    \left\{ x\in C[0,T)\setminus\{0\} : \vert A_s(x)-\alpha\vert\leq \frac{1}{n}\right\}
      \in\cB(C[0,T)),
    \qquad s\in(0,T),\quad n\in\NN,
 \]
 and hence \ $S_\alpha\in\cB(C[0,T))$.
\ For all \ $\alpha, \beta\in\RR$, \ $\alpha\ne\beta$, \ we have
 \ $S_\alpha\cap S_\beta=\emptyset$ \ and, by \eqref{STRONG_CONSISTENCY},
 \begin{align*}
   &\PP_{X^{(\alpha)},\,T}^T(S_\alpha)=\PP(\lim_{t\uparrow T}\halpha_t^{(X^{(\alpha)})}=\alpha)=1,
    \qquad
     \PP_{X^{(\beta)},\,T}^T(S_\beta)=\PP(\lim_{t\uparrow T}\halpha_t^{(X^{(\beta)})}=\beta)=1,\\
   &\PP_{X^{(\alpha)},\,T}^T(S_\beta)=\PP_{X^{(\beta)},\,T}^T(S_\alpha)=0,
 \end{align*}
 which implies the assertion by definition of singularity.
\proofend

In what follows we will study the connections between the technique of the proof of our
 Theorem \ref{THM_singular2} and the very general results on singularity and absolute continuity
 due to Jacod and Shiryaev \cite[Chapter IV]{JacShi}.
In fact, we also present a second proof of Theorem \ref{THM_singular2}.

First we recall that the proof of Theorem \ref{THM_singular2} is based on the strong
 consistency of the MLE of \ $\alpha$, \ see Barczy and Pap \cite[Theorem 3.12]{BarPap2}.
A short outline of the proof of Theorem 3.12 in Barczy and Pap \cite{BarPap2} (specifying for
 \ $\alpha$-Wiener bridges) sounds as follows.
Using the explicit form of the Laplace transform of
 \ $\int_0^t\frac{(X^{(\alpha)}_u)^2}{(T-u)^2}\,\dd u$, \ $t\in[0,T)$, \ due to Barczy and Pap
 \cite[Theorem 4.1]{BarPap2}, one can check that
 \[
  \lim_{t\uparrow T}
      \EE\exp\left\{-\int_0^t\frac{(X_u^{(\alpha)})^2}{(T-u)^2} \,\dd u
               \right\}
   =0,\qquad  \forall\;\alpha\in\RR.
 \]
 Hence
 \begin{align}\label{SEGED_STRONG_CONSISTENCY}
    \PP\left(\lim_{t\uparrow T}\int_0^t\frac{(X^{(\alpha)}_u)^2}{(T-u)^2}\,\dd u=\infty\right)
            =1,\qquad \forall\;\alpha\in\RR,
 \end{align}
 which easily implies strong consistency of the MLE of \ $\alpha$.
\ It will turn out that if we apply Theorem 4.23 in Jacod and Shiryaev \cite[Chapter IV]{JacShi}
 for proving \ $\PP_{X^{(\alpha)},\,T}^T \perp \PP_{X^{(\beta)},\,T}^T$ \ with
 \ $\alpha,\beta\in\RR$, \ $\alpha\ne\beta$, \ then we have to check condition \eqref{SEGED_STRONG_CONSISTENCY}.
We also note that the fact that condition \eqref{SEGED_STRONG_CONSISTENCY} has to be checked
 is in accordance with part (i) of Theorem 1 in Ben-Ari and Pinsky \cite{BenPin}.
But we emphasize that Ben-Ari and Pinsky's result is valid for time-homogeneous diffusions
 and hence we can not use it for $\alpha$-Wiener bridges.
By giving a second proof of Theorem \ref{THM_singular2}, we shed more light on the role of
 condition \eqref{SEGED_STRONG_CONSISTENCY}.

\noindent {\bf Second proof of Theorem \ref{THM_singular2}.}
Let \ $\alpha,\beta\in\RR$, \ $\alpha\ne\beta$ \ be fixed.
Let us introduce the process \ $(\widetilde X^{(\alpha)}_t)_{t\geq 0}$ \ given by
 \[
    \widetilde X^{(\alpha)}_t:= X^{(\alpha)}_{\Psi(t)}, \qquad t\geq 0,
 \]
 where \ $\Psi:[0,\infty)\mapsto[0,T)$, \ $\Psi(t)=\frac{2T}{\pi}\arctan(t)$, \ $t\geq 0$.
\ Then, by the SDE \eqref{alpha_W_bridge} and a change of variable, we get for all \ $t\geq 0$,
 \begin{align*}
   \widetilde X^{(\alpha)}_t
     =-\alpha\int_0^{\Psi(t)}\frac{X^{(\alpha)}_s}{T-s}\,\dd s
        +B_{\Psi(t)}
     =-\alpha\int_0^t\frac{X^{(\alpha)}_{\Psi(u)}}{T-\Psi(u)}\psi(u)\,\dd u
        +B_{\Psi(t)},
 \end{align*}
 where \ $\psi(t):=\frac{\dd}{\dd t}\Psi(t)$, \ $t\geq 0$, \ and
 \ $(B_{\Psi(t)})_{t\geq 0}$ \ is a Wiener process with variance function
 \ $\Psi(t)$, \ $t\geq 0$, \ see Definitions 4.9 in Chapter I in Jacod and Shiryaev \cite{JacShi}.

\noindent Let us consider the filtered space \ $(C[0,\infty),\cB,(\cB_t)_{t\geq 0})$,
 \ where \ $\cB$ \ is the Borel $\sigma$-algebra \ $\cB(C[0,\infty))$ \ on \ $C[0,\infty)$ \ and
 \ $\cB_t$, \ $t\geq 0$, \  defined as follows.
For all \ $t\geq 0$, \ let
 \[
  \cB_t:=\bigcap_{\varepsilon>0}\rho_{t+\varepsilon}^{-1}(\cB),
 \]
 where \ $\rho_t:C[0,\infty)\to C[0,\infty)$ \ defined by \ $(\rho_tx)(s):=x(t\wedge s)$
 \ for \ $s\geq 0$, \ $x\in C[0,\infty)$.
\ Then the filtration \ $(\cB_t)_{t\geq 0}$ \ is right-continuous, since for all \ $t\geq 0$,
 \[
  \cB_t=\bigcap_{\varepsilon>0}\rho_{t+\varepsilon}^{-1}(\cB)
       =\bigcap_{\varepsilon>0}\bigcap_{\delta>0}\rho_{t+\varepsilon+\delta}^{-1}(\cB)
       =\bigcap_{\varepsilon>0}\cB_{t+\varepsilon}.
 \]
Moreover, since \ $\rho_t^{-1}(\cB)\subset\cB_t$ \ for all \ $t\geq 0$, \ by Problem 2.4.2 in
 Karatzas and Shreve \cite{KarShr}, we get
 \begin{align*}
  \cB=\sigma\left(\bigcup_{t\geq0}\cB_t\right).
 \end{align*}
Let \ $\PP_{\widetilde X^{(\alpha)}}$ \ and \ $\PP_{\widetilde X^{(\beta)}}$ \ denote the
 law of the processes \ $(\widetilde X_t^{(\alpha)})_{t\geq 0}$ \ and
 \ $(\widetilde X_t^{(\beta)})_{t\geq 0}$ \ on \ $(C[0,\infty),\cB)$, \ respectively.
\ We check that
 \[
  \PP_{X^{(\alpha)},\,T}^T \perp \PP_{X^{(\beta)},\,T}^T
   \qquad \Longleftrightarrow\qquad
   \PP_{\widetilde X^{(\alpha)}} \perp \PP_{\widetilde X^{(\beta)}}.
 \]
Indeed, by definition, \ $\PP_{X^{(\alpha)},\,T}^T \perp \PP_{X^{(\beta)},\,T}^T$
 \ means that there exist so called distinguishing sets \ $S_\alpha$ \ and
 \ $S_\beta$ \ in \ $\cB(C[0,T))$ \ such that \ $S_\alpha\cap S_\beta=\emptyset$ \ and
 \[
   \PP\big(\{\omega\in\Omega : (X^{(\alpha)}_t(\omega))_{t\in[0,T)}\in S_\alpha\}\big)
       =\PP\big(\{\omega\in\Omega : (X^{(\beta)}_t(\omega))_{t\in[0,T)}\in S_\beta\}\big)
       =1.
 \]
Similarly, \ $\PP_{\widetilde X^{(\alpha)}} \perp \PP_{\widetilde X^{(\beta)}}$
 \ means that there exist \ $\widetilde S_\alpha$ \ and \ $\widetilde S_\beta$ \ in
 \ $\cB$ \ such that \ $\widetilde S_\alpha\cap \widetilde S_\beta=\emptyset$ \ and
 \[
   \PP\big(\{\omega\in\Omega : (\widetilde X^{(\alpha)}_t(\omega))_{t\geq 0}\in \widetilde S_\alpha\}\big)
     =\PP\big(\{\omega\in\Omega : (\widetilde X^{(\beta)}_t(\omega))_{t\geq 0}\in \widetilde S_\beta\}\big)
     =1.
 \]
Using that
 \begin{align*}
   \{\omega\in\Omega : (\widetilde X^{(\alpha)}_t(\omega))_{t\geq 0}\in \widetilde S_\alpha\}
     &=\{\omega\in\Omega : (X^{(\alpha)}_{\Psi(t)}(\omega))_{t\geq 0}\in \widetilde S_\alpha\}\\
     &=\{\omega\in\Omega : (X^{(\alpha)}_t(\omega))_{t\in[0,T)}\in \Psi^{-1}(\widetilde S_\alpha)\},
 \end{align*}
 where \ $\Psi^{-1}(\widetilde S_\alpha):=\{f\circ\Psi^{-1} : f\in\widetilde S_\alpha\}$,
 \ singularity of \ $\PP_{\widetilde X^{(\alpha)}}$ \ and \ $\PP_{\widetilde X^{(\beta)}}$
 \ with distinguishing sets \ $\widetilde S_\alpha, \widetilde S_\beta\in\cB$ \ implies
 singularity of \ $\PP_{X^{(\alpha)},\,T}^T$ \ and \ $\PP_{X^{(\beta)},\,T}^T$
 \ with distinguishing sets \ $\Psi^{-1}(\widetilde S_\alpha), \Psi^{-1}(\widetilde S_\beta)\in\cB(C[0,T))$.
\ The converse statement can be thought over similarly.

\noindent Hence by Corollary 2.8 in Chapter IV in Jacod and Shiryaev \cite{JacShi},
 to prove the assertion it is enough to check that the measures \ $\PP_{\widetilde X^{(\alpha)}}$
 \ and \ $\PP_{\widetilde X^{(\beta)}}$ \ are locally equivalent with respect to each other
 (where the restrictions of the measures refers to the given filtration
 \ $(\cB_t)_{t\geq 0}$) \ and that \ $\PP_{\widetilde X^{(\alpha)}}(\lim_{t\to\infty}h_t^{(1/2)}<\infty)=0$,
 \ where \ $(h_t^{(1/2)})_{t>0}$ \ is the Hellinger process of order \ $1/2$
 \ between \ $\PP_{\widetilde X^{(\alpha)}}$ \ and \ $\PP_{\widetilde X^{(\beta)}}$.
\ Using that the continuity of the process \ $X^{(\alpha)}$ \ implies that
 the process
 \[
   \int_0^t \frac{(X^{(\alpha)}_{\Psi(u)})^2}{(T-\Psi(u))^2}\psi(u)^2 \,\dd u,
     \qquad t\geq 0,
 \]
 does not jump to infinity (for the definition of jumping to infinity, see, e.g., Definitions 5.8
 (ii) in Chapter III in Jacod and Shiryaev \cite{JacShi}), by \eqref{SEGED13} and
 a generalization of part (b) and (c) of Theorem 4.23 in Chapter IV in Jacod and Shiryaev \cite{JacShi},
 we have the measures \ $\PP_{\widetilde X^{(\alpha)}}$ \ and
 \ $\PP_{\widetilde X^{(\beta)}}$ \ are locally equivalent with respect to each other
 and the process
 \begin{align}\label{SEGED_JACOD_SHIRYAEV}
     \frac{(\alpha-\beta)^2}{8}\int_0^t\frac{x(\Psi(u))^2}
                               {(T-\Psi(u))^2}\psi(u)\,\dd u,
       \qquad x\in C[0,\infty),\quad t>0,
 \end{align}
 is a version of the Hellinger process \ $(h_t^{(1/2)})_{t>0}$.
\ Indeed, using the notations of Sections 3a and 4b in Chapter IV in Jacod and Shiryaev \cite{JacShi},
 we have \ $C(t)=\Psi(t)=\int_0^t\psi(s)\,\dd s$, \ $t\geq 0$, \ and
 \begin{align*}
  &\beta_s(x)=-\alpha\frac{x(\Psi(s))}{T-\Psi(s)}\psi(s),
         \qquad x\in C[0,\infty),\quad s\geq 0,\\
  &\beta_s'(x)=-\beta\frac{x(\Psi(s))}{T-\Psi(s)}\psi(s),
         \qquad x\in C[0,\infty),\quad s\geq 0,\\
  &\widetilde\beta_s(s)
           =\frac{\beta_s(x)-\beta_s'(x)}{\psi(s)}
           =-(\alpha-\beta)\frac{x(\Psi(s))}{T-\Psi(s)},
         \qquad x\in C[0,\infty),\quad s\geq 0.
 \end{align*}
Hence using the very same arguments given in the proof of Theorem 4.23 in Jacod and Shiryaev
 \cite[Chapter IV]{JacShi}, we get \eqref{SEGED_JACOD_SHIRYAEV}.
By a change of variable, we have
 \[
  \int_0^t\frac{x(\Psi(u))^2}{(T-\Psi(u))^2}\psi(u)\,\dd u
    =\int_0^{\Psi(t)}\frac{x(s)^2}{(T-s)^2}\,\dd s,
            \qquad x\in C[0,\infty),\quad t\geq 0,
 \]
 and hence \ $\PP_{\widetilde X^{(\alpha)}}(\lim_{t\to\infty}h_t^{(1/2)}<\infty)=0$, \ $\alpha\in\RR$,
 \ is equivalent with \eqref{SEGED_STRONG_CONSISTENCY}, i.e., to prove the assertion it is enough to
 verify \eqref{SEGED_STRONG_CONSISTENCY}.
As it was mentioned earlier, as a special case of the proof of Theorem 3.12 in Barczy and Pap
 \cite{BarPap2} we get \eqref{SEGED_STRONG_CONSISTENCY}.
\proofend

\begin{Rem}
If \ $\alpha, \beta\in\RR$, \ $\alpha\ne\beta$, \ by Theorem \ref{THM_singular2}, we have
 \ $\PP^T_{X^{(\alpha)},\,T} \perp \PP^T_{X^{(\beta)},\,T}$.
\ Moreover, in the proof of the theorem, we also constructed disjoint sets \ $S_\alpha$ \ and \ $S_\beta$
 \ in \ $\cB(C[0,T))$ \ that distinguish between the measures
 in the sense that \ $\PP^T_{X^{(\alpha)},\,T}(S_\alpha)=1$ \ and  \ $\PP^T_{X^{(\beta)},\,T}(S_\beta)=1$.
\ We note that for some special time-homogeneous ($1$-dimensional) diffusions Ben-Ari and Pinsky
 \cite[Propositions 1, 2 and 3]{BenPin} also gave "illuminating" distinguishing sets.
\proofend
\end{Rem}

\begin{Rem}\label{REMARK3}
In case of \ $\alpha\geq0$, \ by Lemma \ref{LEMMA5}, one can define a probability measure
 \ $\PP_{Y^{(\alpha)},\,T}^T$ \ on \ $(C[0,T],\cB(C[0,T]))$, \ as the law of the process
 \ $(Y^{(\alpha)}_t)_{t\in[0,T]}$ \ given in Lemma \ref{LEMMA5}.
Then, by Theorem \ref{THM_singular2}, for all \ $\alpha, \beta\geq 0$,
 \ $\alpha\ne\beta$, \ we get \ $\PP_{Y^{(\alpha)},\,T}^T \perp \PP_{Y^{(\beta)},\,T}^T$.
\ Indeed, since \ $\PP_{X^{(\alpha)},\,T}^T \perp \PP_{X^{(\beta)},\,T}^T$,
 \ there exist sets \ $S_\alpha$ \ and \ $S_\beta$ \ in \ $\cB(C[0,T))$ \ such that
 \ $S_\alpha\cap S_\beta=\emptyset$ \ and
 \ $\PP_{X^{(\alpha)},\,T}^T(S_\alpha)=\PP_{X^{(\beta)},\,T}^T(S_\beta)=1$.
\ For all \ $B\in\cB(C[0,T))$, \ let us introduce the notation
 \[
   \widetilde B
     :=\Big\{ x\in C[0,T] : x\vert_{[0,T)}\in B \quad\text{and}\quad
                \exists \, \lim_{t\uparrow T}x(t)\in\RR \Big\}.
 \]
Then we have \ $\widetilde S_\alpha, \widetilde S_\beta\in\cB(C[0,T])$,
 \ $\widetilde S_\alpha\cap \widetilde S_\beta=\emptyset$ \ and
 \ $\PP_{Y^{(\alpha)},\,T}^T(\widetilde S_\alpha)=\PP_{Y^{(\beta)},\,T}^T(\widetilde S_\beta)=1$.
\ As a special case, we also have for all \ $\alpha>0$, \ the probability measure
 \ $\PP_{Y^{(\alpha)},\,T}^T$ \ and the standard Wiener measure \ $\PP_{Y^{(0)},\,T}$
 \ on \ $(C[0,T],\cB(C[0,T]))$ \ are singular.
\proofend
\end{Rem}

\begin{Rem}
We note that Theorem \ref{THM_singular2} is not an astonishing result.
One can easily formulate conditions on a general time-inhomogeneous diffusion process
 under which the same kind of singularity holds.
Namely, let us consider a process \ $(X^{(\theta)}_t)_{t\geq 0}$ \ given by the SDE
 \begin{align}\label{SDE_general_inhomogeneous}
   \begin{cases}
     \dd X^{(\theta)}_t=\theta a(t,X^{(\theta)}_t)\,\dd t+\dd B_t,\qquad t\geq 0,\\
     \phantom{\dd} X_0^{(\theta)}=0,
   \end{cases}
 \end{align}
 where \ $a:[0,\infty)\times\RR\to\RR$ \ is a known Borel-measurable function,
 \ $(B_t)_{t\geq 0}$ \ is a standard Wiener process, and \ $\theta\in\RR$ \ is an unknown parameter.
Let us suppose that the SDE \eqref{SDE_general_inhomogeneous} has a unique strong solution
 \ $(X^{(\theta)}_t)_{t\geq 0}$ \ for all \ $\theta\in\RR$.
\ For all \ $\theta\in\RR$, \ let us denote by \ $\PP_\theta$ \ the law of \ $(X^{(\theta)}_t)_{t\geq 0}$
 \ on \ $(C[0,\infty),\cB(C[0,\infty)))$.
\ Let us suppose that for all \ $t>0$ \ and all \ $\theta\in\RR$,
 \begin{align*}
    \PP\left(\int_0^t a(s,X^{(\theta)}_s)^2\,\dd s<\infty\right)=1.
 \end{align*}
As it is explained in details in the second proof of Theorem \ref{THM_singular2},
 using Theorem 4.23 in Chapter IV in Jacod and Shiryaev \cite{JacShi}, we get
 for all \ $\theta_1$, $\theta_2\in\RR$, \ $\theta_1\ne\theta_2$,
 \[
  \PP_{\theta_1} \perp \PP_{\theta_2}
    \qquad  \Longleftrightarrow\qquad
   \PP\left(\lim_{t\to \infty}
      \int_0^t a(u,X^{(\theta_i)}_u)^2 \,\dd u
             =\infty\right)=1
      \qquad \text{for \ $i=1$ \ or \ $i=2$.}
 \]

\noindent Concerning singularity of \ $\PP_{\theta_1}$ \ and \ $\PP_{\theta_2}$,
 \ the point is that whether the imposed conditions can be checked for a given diffusion process.
And in this respect, time-inhomogeneous diffusions in general represent a hard task.
\proofend
\end{Rem}

Concerning the SDE \eqref{alpha_W_bridge} one can ask why the diffusion coefficient
 in the SDE \eqref{alpha_W_bridge} is identically 1.
The point is only that it is supposed to be a known and positive constant.
Remember that in many cases the measures induced by processes with different diffusion coefficients
 are singular and continuous-time statistical inference for this type of model is often trivial.
In what follows we consider this phenomenon in details.
For all \ $T\in(0,\infty)$, \ $\alpha\in\RR$ \ and \ $\sigma>0$, \ let us introduce
 the time-inhomogeneous diffusion process \ $(X_t^{(\alpha,\sigma)})_{t\in[0,T)}$ \ given by the SDE
 \begin{align}\label{alpha_W_bridge_with_sigma}
  \begin{cases}
   \dd X_t^{(\alpha,\sigma)}=-\frac{\alpha}{T-t}\,X_t^{(\alpha,\sigma)}\,\dd t+\sigma\dd B_t,\qquad t\in[0,T),\\
   \phantom{\dd} X_0^{(\alpha,\sigma)}=0,
  \end{cases}
 \end{align}
 where \ $(B_t)_{t\geq 0}$ \ is a 1-dimensional standard Wiener process.
By formula (5.6.6) in Karatzas and Shreve \cite{KarShr}, the SDE \eqref{alpha_W_bridge_with_sigma}
 has a unique strong solution, namely,
 \[
    X_t^{(\alpha,\sigma)}=
      \sigma\int_0^t\left(\frac{T-t}{T-s}\right)^\alpha\,\dd B_s,\qquad t\in[0,T).
 \]
For all \ $t\in(0,T)$, \ let \ $\PP_{X^{(\alpha,\sigma)},\,t}$ \ be the law of the process
 \ $(X_s^{(\alpha,\sigma)})_{s\in[0,t]}$ \ given by the SDE \eqref{alpha_W_bridge_with_sigma}
 on \ $\big(C[0,t],\cB(C[0,t])\big)$.

\begin{Thm}\label{THM_singular1}
For all \ $\alpha_1, \alpha_2\in\RR$, \ $\sigma_1>0$, \ $\sigma_2>0$ \ and \ $t\in(0,T)$,
 \ the following dichotomy holds:
 \begin{align*}
     &\PP_{X^{(\alpha_1,\sigma_1)},\,t}\sim \PP_{X^{(\alpha_2,\sigma_2)},\,t}
         \qquad\text{if \ $\sigma_1=\sigma_2$, }\\
     &\PP_{X^{(\alpha_1,\sigma_1)},\,t}\perp \PP_{X^{(\alpha_2,\sigma_2)},\,t}
         \qquad\text{if \ $\sigma_1\ne \sigma_2$.}
 \end{align*}
\end{Thm}

\noindent{\bf First proof.}
 In case of \ $\sigma_1=\sigma_2$,  \ the equivalence of \ $\PP_{X^{(\alpha_1,\sigma_1)},\,t}$ \ and
 \ $\PP_{X^{(\alpha_2,\sigma_2)},\,t}$ \ follows from Theorem 7.20 or Theorem 7.19
 in Liptser and Shiryaev \cite{LipShiII} and from \eqref{SEGED13}.

\noindent  Let us suppose now that \ $\sigma_1\ne\sigma_2$. \ For all
 \ $\alpha\in\RR$ \ and \ $\sigma>0$,  \ by giving a direct proof,
 we show the following Baxter type result
 \begin{align}\label{SEGED_ortogonality_1}
   \PP\left(\lim_{n\to\infty}\sum_{j=1}^n(X_{jt/n}^{(\alpha,\sigma)}-X_{(j-1)t/n}^{(\alpha,\sigma)})^2
           =t\sigma^2\right)=1.
 \end{align}
By the SDE \eqref{alpha_W_bridge_with_sigma}, we have
 \begin{align}\label{SEGED_ortogonality_2}
  \begin{split}
   \sum_{j=1}^n(X_{jt/n}^{(\alpha,\sigma)}&-X_{(j-1)t/n}^{(\alpha,\sigma)})^2\\
      =&\alpha^2\sum_{j=1}^n
         \left(\int_0^{jt/n}\frac{X_u^{(\alpha,\sigma)}}{T-u}\,\dd u
                -\int_0^{(j-1)t/n}\frac{X_u^{(\alpha,\sigma)}}{T-u}\,\dd u\right)^2\\
        &-2\alpha\sigma\sum_{j=1}^n
           \left(\int_0^{jt/n}\frac{X_u^{(\alpha,\sigma)}}{T-u}\,\dd u
                -\int_0^{(j-1)t/n}\frac{X_u^{(\alpha,\sigma)}}{T-u}\,\dd u\right)
                   (B_{jt/n}-B_{(j-1)t/n})\\
       &+\sigma^2\sum_{j=1}^n(B_{jt/n}-B_{(j-1)t/n})^2.
  \end{split}
 \end{align}
It is known that
 \begin{align}\label{SEGED_ortogonality_3}
   \PP\left(\lim_{n\to\infty}\sum_{j=1}^n(B_{jt/n}-B_{(j-1)t/n})^2=t\right)=1,
 \end{align}
 see, e.g., Lemma 4.3 in Liptser and Shiryaev \cite{LipShiI}.
Moreover, by Lagrange's mean value theorem, one can think it over that
 for all \ $t\in[0,T)$ \ and for all continuous functions \ $f:[0,t]\to\RR$, \ we have
 \[
   \lim_{n\to\infty}\sum_{j=1}^n
      \left(\int_0^{jt/n}f(x)\,\dd x-\int_0^{(j-1)t/n}f(x)\,\dd x\right)^2
      =0.
 \]
Since for all \ $t\in(0,T)$, \ the process \ $\left(\frac{X_u^{(\alpha,\sigma)}}{T-u}\right)_{u\in[0,t]}$
 \ is continuous, we have
 \begin{align}\label{SEGED_ortogonality_4}
      \PP\left(\lim_{n\to\infty}\sum_{j=1}^n
            \left(\int_0^{jt/n}\frac{X_u^{(\alpha,\sigma)}}{T-u}\,\dd u
                   -\int_0^{(j-1)t/n}\frac{X_u^{(\alpha,\sigma)}}{T-u}\,\dd u\right)^2=0\right)=1.
 \end{align}
Now we check that
 \begin{align}\label{SEGED_ortogonality_5}
     \PP\left(\lim_{n\to\infty}\sum_{j=1}^n
                 \left(\int_0^{jt/n}\frac{X_u^{(\alpha,\sigma)}}{T-u}\,\dd u
                   -\int_0^{(j-1)t/n}\frac{X_u^{(\alpha,\sigma)}}{T-u}\,\dd u\right)
                    (B_{jt/n}-B_{(j-1)t/n})=0\right)=1.
 \end{align}
By Cauchy--Schwartz's inequality, we have
 \begin{align*}
   \sum_{j=1}^n
       &\left(\int_0^{jt/n}\frac{X_u^{(\alpha,\sigma)}}{T-u}\,\dd u
                   -\int_0^{(j-1)t/n}\frac{X_u^{(\alpha,\sigma)}}{T-u}\,\dd u\right)
        (B_{jt/n}-B_{(j-1)t/n})\\
       &\leq\sqrt{\sum_{j=1}^n
         \left(\int_0^{jt/n}\frac{X_u^{(\alpha,\sigma)}}{T-u}\,\dd u
                   -\int_0^{(j-1)t/n}\frac{X_u^{(\alpha,\sigma)}}{T-u}\,\dd u\right)^2}
         \sqrt{\sum_{j=1}^n(B_{jt/n}-B_{(j-1)t/n})^2},
 \end{align*}
 with probability one, and then \eqref{SEGED_ortogonality_3}
 and \eqref{SEGED_ortogonality_4} implies \eqref{SEGED_ortogonality_5}.
By \eqref{SEGED_ortogonality_3}, \eqref{SEGED_ortogonality_4} and \eqref{SEGED_ortogonality_5},
 using \eqref{SEGED_ortogonality_2} we have \eqref{SEGED_ortogonality_1}.
Then, using the definition of singularity of measures, \eqref{SEGED_ortogonality_1} implies that
 \ $\PP_{X^{(\alpha,\sigma_1)},\,t} \perp \PP_{X^{(\alpha,\sigma_2)},\,t}$ \ for all
 \ $\alpha\in\RR$ \ and \ $\sigma_1\ne\sigma_2$.
\ In case of \ $\alpha_1\ne\alpha_2$ \ and \ $\sigma_1\ne\sigma_2$ \ we have
 \ $\PP_{X^{(\alpha_1,\sigma_1)},\,t} \sim \PP_{X^{(\alpha_2,\sigma_1)},\,t}$ \
 and \ $\PP_{X^{(\alpha_2,\sigma_1)},\,t} \perp \PP_{X^{(\alpha_2,\sigma_2)},\,t}$,
 \ which imply that \ $\PP_{X^{(\alpha_1,\sigma_1)},\,t} \perp \PP_{X^{(\alpha_2,\sigma_2)},\,t}$.

\noindent{\bf Second proof.}
Using Baxter's theorem due to Baxter \cite[Theorem 1]{Bax}, we show that
 for all \ $\alpha\in\RR$ \ and \ $\sigma>0$,
 \begin{align}\label{SEGED_BAXTER}
    \PP\left(\lim_{n\to\infty}\sum_{k=1}^{2^n}
            \left(X^{(\alpha,\sigma)}_{kt/2^n}-X^{(\alpha,\sigma)}_{(k-1)t/2^n}\right)^2
                     =t\sigma^2\right)=1,
 \end{align}
 which is also enough (like \eqref{SEGED_ortogonality_1}) to ensure that
 \ $\PP_{X^{(\alpha,\sigma_1)},\,t} \perp \PP_{X^{(\alpha,\sigma_2)},\,t}$ \ for all
 \ $\alpha\in\RR$ \ and \ $\sigma_1\ne\sigma_2$.
\ For all \ $t\in(0,T)$, \ $(X^{(\alpha,\sigma)}_s)_{s\in[0,t]}$ \ is a Gauss process
 with identically 0 mean function, and to have right to apply Baxter's theorem,
 we need to check that the covariance function of \ $(X^{(\alpha,\sigma)}_s)_{s\in[0,t]}$ \ is continuous
 on \ $[0,t]\times[0,t]$ \ and has uniformly bounded second derivatives on
 \ $[0,t]\times[0,t]\setminus\{(s,s):s\in[0,t]\}$.

\noindent In case of \ $\alpha\ne\frac{1}{2}$,  \ by \eqref{SEGED_COV2}, we get for all \ $u,\,v\in[0,t]$,
 \begin{align*}
   r(v,u):=\cov(X^{(\alpha,\sigma)}_v,X^{(\alpha,\sigma)}_u)
          =\sigma^2
           \frac{(T-v)^{\alpha}(T-u)^{\alpha}}{1-2\alpha}
           \Big(T^{1-2\alpha}-(T-(v\wedge u))^{1-2\alpha}\Big).
 \end{align*}
Then \ $r(v,u)$, \ $v,u\in[0,t]$, \ is continuous, and clearly, if \ $0\leq v<u\leq t$, \ then
 \begin{align*}
   &\frac{\partial r}{\partial v}(v,u)
      =\sigma^2\frac{T^{1-2\alpha}}{1-2\alpha}(T-u)^\alpha\alpha(T-v)^{\alpha-1}(-1)
        -\frac{\sigma^2}{1-2\alpha}(T-u)^{\alpha}(1-\alpha)(T-v)^{-\alpha}(-1),\\[1mm]
   &\frac{\partial r}{\partial u}(v,u)
      =\sigma^2\frac{T^{1-2\alpha}}{1-2\alpha}(T-v)^\alpha\alpha(T-u)^{\alpha-1}(-1)
        -\frac{\sigma^2}{1-2\alpha}(T-v)^{1-\alpha}\alpha(T-u)^{\alpha-1}(-1).
 \end{align*}
If \ $0\leq u<v\leq t$, \ then
 \begin{align*}
   &\frac{\partial r}{\partial v}(v,u)
      =\sigma^2\frac{T^{1-2\alpha}}{1-2\alpha}(T-u)^\alpha\alpha(T-v)^{\alpha-1}(-1)
        -\frac{\sigma^2}{1-2\alpha}(T-u)^{1-\alpha}\alpha(T-v)^{\alpha-1}(-1),\\[1mm]
   &\frac{\partial r}{\partial u}(v,u)
       =\sigma^2\frac{T^{1-2\alpha}}{1-2\alpha}(T-v)^\alpha\alpha(T-u)^{\alpha-1}(-1)
         -\frac{\sigma^2}{1-2\alpha}(T-v)^{\alpha}(1-\alpha)(T-u)^{-\alpha}(-1).
 \end{align*}
This implies that \ $(X^{(\alpha,\sigma)}_s)_{s\in[0,t]}$ \ has uniformly bounded first
 derivatives on \ $[0,t]\times[0,t]\setminus\{{(s,s):s\in[0,t]}\}$, \ and similarly
 one can check that the second derivatives also admit this property.
\ Moreover, for all \ $u\in[0,t)$,
 \begin{align*}
    D^+(u):&=\lim_{v\downarrow u}\frac{r(u,u)-r(v,u)}{u-v}
            =\lim_{v\downarrow u}\frac{r(v,u)-r(u,u)}{v-u}
            =\lim_{v\downarrow u}\frac{\partial r}{\partial v}(v,u)\\
           &=-\frac{\sigma^2\alpha}{1-2\alpha}T^{1-2\alpha}(T-u)^{2\alpha-1}
              +\frac{\sigma^2\alpha}{1-2\alpha}.
  \end{align*}
Similarly, for all \ $u\in(0,t]$,
 \begin{align*}
    D^-(u):&=\lim_{v\uparrow u}\frac{r(u,u)-r(v,u)}{u-v}
            =\lim_{v\uparrow u}\frac{r(v,u)-r(u,u)}{v-u}
            =\lim_{v\uparrow u}\frac{\partial r}{\partial v}(v,u)\\
           &=-\frac{\sigma^2\alpha}{1-2\alpha}T^{1-2\alpha}(T-u)^{2\alpha-1}
              +\frac{\sigma^2(1-\alpha)}{1-2\alpha}.
 \end{align*}
Hence for all \ $t\in[0,T)$,
 \[
      \int_0^t(D^-(u)-D^+(u))\,\dd u
         =\sigma^2\int_0^t\left(\frac{1-\alpha}{1-2\alpha}-\frac{\alpha}{1-2\alpha}\right)\,\dd u
         =\sigma^2\int_0^t1\,\dd u
         =\sigma^2 t,
 \]
 and then Baxter's theorem due to Baxter \cite[Theorem 1]{Bax} yields \eqref{SEGED_BAXTER}.

In case of \ $\alpha=\frac{1}{2}$, \ similarly to the case \ $\alpha\ne\frac{1}{2}$,
 \ one can check that the covariance function of \ $(X^{(1/2,\sigma)}_s)_{s\in[0,t]}$
 \ is continuous on \ $[0,t]\times[0,t]$, \ and has uniformly bounded first and second
 derivatives on \ $[0,t]\times[0,t]\setminus\{(s,s):s\in[0,t]\}$.
\ Moreover, by \eqref{SEGED_COV2}, for all \ $u\in[0,t)$,
 \begin{align*}
    D^+(u)
       &=\lim_{v\downarrow u}
          \frac{\sigma^2\sqrt{(T-v)(T-u)}\ln\left(\frac{T}{T-u}\right)
                 - \sigma^2(T-u)\ln\left(\frac{T}{T-u}\right)}
               {v-u}\\
       &=\sigma^2\ln\left(\frac{T}{T-u}\right)\sqrt{T-u}
          \,\frac{\dd}{\dd u}\sqrt{T-u}
       =-\frac{\sigma^2}{2}\ln\left(\frac{T}{T-u}\right).
 \end{align*}
Similarly, by \eqref{SEGED_COV2}, for all \ $u\in(0,t]$,
 \begin{align*}
    D^-(u)
       &=\lim_{v\uparrow u}
          \frac{\sigma^2\sqrt{(T-v)(T-u)}\ln\left(\frac{T}{T-v}\right)
                 - \sigma^2(T-u)\ln\left(\frac{T}{T-u}\right)}
               {v-u}\\
       &=\sigma^2\sqrt{T-u}
          \frac{\dd}{\dd u}\left(\sqrt{T-u}\ln\left(\frac{T}{T-u}\right)\right)
       =-\frac{\sigma^2}{2}\ln\left(\frac{T}{T-u}\right)+\sigma^2.
 \end{align*}
Hence for all \ $t\in[0,T)$,
 \[
      \int_0^t(D^-(u)-D^+(u))\,\dd u
         =\sigma^2\int_0^t1\,\dd u
         =\sigma^2 t,
 \]
 and then Baxter's theorem due to Baxter \cite[Theorem 1]{Bax} yields \eqref{SEGED_BAXTER}.
\proofend

We note that the same dichotomy that we have in Theorem \ref{THM_singular1} holds for Ornstein--Uhlenbeck
 processes, see, e.g., page 226 in Arat\'o, Pap and van Zuijlen \cite{AraPapZui}.

\begin{Rem}
For all \ $\alpha\in\RR$ \ and \ $\sigma>0$, \ let \ $\PP_{X^{(\alpha,\sigma)}}$
 \ denote the law of the process \ $(X^{(\alpha,\sigma)}_s)_{s\in[0,T)}$ \ on
 \ $(C[0,T),\cB(C[0,T)))$.
\ By the proof of Theorem \ref{THM_singular2}, we get
 \ $\PP_{X^{(\alpha_1,\sigma_1)}} \perp \PP_{X^{(\alpha_2,\sigma_2)}}$
 \ for all \ $\alpha_1,\alpha_2\in\RR$ \ and \ $\sigma_1>0$, $\sigma_2>0$.
\proofend
\end{Rem}

\vskip0.7cm

\noindent{\bf\large Acknowledgements.}
The first author has been supported by the Hungarian Scientific Research Fund under
 Grants No.\ OTKA--F046061/2004 and OTKA T-048544/2005.
The second author has been supported by the Hungarian Scientific Research Fund under
 Grant No.\ OTKA T-048544/2005.

\end{document}